\documentclass{article}

\usepackage{amssymb}
\usepackage{amsmath}
\usepackage[mathscr]{eucal}
\usepackage{epsfig}

\newtheorem{Lemma}{Lemma}
\newtheorem{Proposition}[Lemma]{Proposition}
\newtheorem{Theorem}[Lemma]{Theorem}

\newtheorem{Fact}{Fact}
\newtheorem{Definition}{Definition}

\newcommand{\cd}{\ \stackrel{d}{\longrightarrow} \ }

\newcommand{\ed}{\ \stackrel{d}{=} \ }

\newcommand{\PP}{\mbox{${\mathcal P}$}}

\newcommand{\EE}{\mbox{${\mathcal E}$}}
\newcommand{\VV}{\mbox{${\mathcal V}$}}
\newcommand{\GG}{\mbox{${\mathcal G}$}}

\newcommand{\Nbold}{\mbox{${\mathbb N}$}}
\newcommand{\Rbold}{\mbox{${\mathbb R}$}}

\newcommand{\Tbold}{\mbox{${\mathbb T}$}}

\newcommand{\bE}{{\bf E}}
\newcommand{\bP}{{\bf P}}

\newcommand{\bi}{{\bf i}}

\newcommand{\sfrac}[2]{{\textstyle\frac{#1}{#2}}}

\newcommand{\qed}{\hfill \ \ \rule{1ex}{1ex}}

\newcommand{\Hbar}{{\overline H}}

\newcommand{\mudiag}{\mu^{\nearrow}}

\newcommand{\stleq}{\preccurlyeq}

\newcommand{\lleq}{\leq^*}

\title{{\bf Bivariate Uniqueness and Endogeny for the Logistic Recursive
            Distributional Equation}}

\author{{\bf Antar Bandyopadhyay} \\
\\
\\
Deaprtment of Mathematics \\
Chalmers University of Technology \\
SE - 412 96, G\"{o}teborg \\ 
SWEDEN \\
E-Mail : {\tt antar@math.chalmers.se}
}

\begin{document}

\maketitle
\bibliographystyle{plain} 

\begin{abstract}

In this article we prove the \emph{bivariate uniqueness} property for
a particular ``max-type'' \emph{recursive distributional equation} (RDE).
Using the general theory developed in \cite{AlBa05} we then show that 
the corresponding \emph{recursive tree process} (RTP) has no 
\emph{external} randomness, more preciously, the RTP is 
\emph{endogenous}. The RDE we consider is so called the 
\emph{Logistic RDE}, which appears in the proof of the $\zeta(2)$-limit 
of the random assignment problem \cite{Al01} using the local weak 
convergence method. Thus this work provides a non-trivial application of
the general theory developed in \cite{AlBa05}. 

\end{abstract}

\vspace{0.2in}


\vspace{0.1in}


\newpage

\section{Introduction and the Main Result}
\label{Sec:Intro}
Fixed-point equations or distributional identities have appeared in the
probability literature for quite a long time in a variety of settings. 
The recent survey of Aldous and Bandyopadhyay \cite{AlBa05}
provides a general framework to study certain type of distributional 
equations. 

Given a space $S$ write $\PP\left(S\right)$ for the set of all probabilities
on $S$. A \emph{recursive distributional equation} (RDE) \cite{AlBa05} is a 
fixed-point equation on $\PP\left(S\right)$ defined as
\begin{equation}
X \ed g \left( \xi ; \left( X_j : 1 \leq j \lleq N \right) \right)
\,\,\,\, \mbox{on} \,\,\, S, 
\label{Equ:RDE}
\end{equation}
where it is assumed that
$\left(X_j\right)_{j \geq 1}$ are i.i.d. $S$-valued random variables
with same distribution as
$X$, and are independent of the pair $\left( \xi, N \right)$. Here
$N$ is a non-negative integer valued random variable, 
which may take the value $\infty$, and $g$ is a given $S$-valued function. 
(In the above equation by ``$\lleq N$'' we mean the left hand side is 
``$ \leq N$'' if $N < \infty$, and ``$< N$'' otherwise). 
In (\ref{Equ:RDE}) the distribution of $X$ is \emph{unknown}, while
the distribution of the pair $\left( \xi, N \right)$ and the function
$g$ are the \emph{known} quantities. Perhaps a more conventional 
(analytic) way of writing the equation (\ref{Equ:RDE}) would be
\begin{equation}
\mu = T\left(\mu\right) \, , 
\label{Equ:RDE-T}
\end{equation}
where $T : \PP \rightarrow \PP\left(S\right)$ is a function defined on
$\PP \subseteq \PP\left(S\right)$ such that $T\left(\mu\right)$ is
the distribution of the right-hand side of the equation (\ref{Equ:RDE}), 
when $\left(X_j\right)_{j \geq 1}$ are i.i.d. $\mu \in \PP$. 

As outlined in \cite{AlBa05} 
in many applications RDEs play a very crucial 
role. Examples include study of Galton-Watson branching processes and
related random trees, probabilistic analysis of algorithms with
suitable recursive structure \cite{Ros92, FiJa00, RosRu01}, 
statistical physics models on trees 
\cite{Al00, AlSt04, Ga04, Ba05a, Ba05b, BaGa05},
and statistical physics and algorithmic questions in the mean-field model
of distance \cite{Al92a, Al01, AlSt04}.
In many of these applications, particularly in the last 
two types mentioned above, often one needs to construct 
a particular tree indexed \emph{stationary}
process related to a given RDE, which is called a
\emph{recursive tree process} (RTP) \cite{AlBa05}.
More precisely, suppose the RDE (\ref{Equ:RDE}) has a 
solution, say $\mu$. Then as shown in \cite{AlBa05}, using the
consistency theorem of Kolmogorov \cite{Bill95},
one can construct a process, say 
$\left(X_{\bi}\right)_{\bi \in \VV}$, indexed by 
$\VV := 
\left(\cup_{d \geq 1} \Nbold^d \right) \cup \left\{ \emptyset \right\}$, 
such that 
\begin{equation}
\begin{array}{cl}
\mbox{(i)} & X_{\bi} \sim \mu \,\,\,\, \forall \,\,\, \bi \in \VV, \\
\mbox{(ii)} & 
\mbox{For each\ \ } d \geq 0, \left(X_{\bi}\right)_{\vert \bi \vert = d}
\mbox{\ \ are independent}, \\
\mbox{(iii)} & 
X_{\bi} = g \left( \xi_{\bi} ; \left( X_{\bi j} : 
1 \leq j \lleq N_{\bi} \right) \right) \,\,\,\, \forall \,\,\, \bi \in \VV, \\
\mbox{(iv)} &
X_{\bi} \mbox{\ \ is independent of\ \ } 
\left\{ \left(\xi_{\bi'}, N_{\bi'} \right) \,\Big\vert\, 
\vert \bi' \vert < \vert \bi \vert \, \right\}
\,\,\,\, \forall \,\,\, \bi \in \VV,
\end{array}
\label{Equ:RTP}
\end{equation}
where $\left(\xi_{\bi}, N_{\bi}\right)_{\bi \in \VV}$ are taken to be i.i.d. 
copies of the pair $\left( \xi, N \right)$, and by 
$\vert \cdot \vert$ we mean the length of a finite word. 
The process $\left(X_{\bi}\right)_{\bi \in \VV}$ is called an invariant 
\emph{recursive tree process} (RTP) with marginal $\mu$. The i.i.d.
random variables $\left(\xi_{\bi}, N_{\bi} \right)_{\bi \in \VV}$ are 
called the \emph{innovation process}. In some sense an invariant 
RTP with marginal $\mu$, 
is an almost sure representation of a solution $\mu$ of the RDE
(\ref{Equ:RDE}). 
Here we note that there is a natural tree structure on $\VV$. 
Taking $\VV$ as the vertex set, we join two words $\bi, \bi' \in \VV$
by an edge, if and only if, $\bi' = \bi j$
or $\bi = \bi' j$, for some $j \in \Nbold$. We will denote this
tree by $\Tbold_{\infty}$. 
The empty-word $\emptyset$ will
be taken as the root of the tree $\Tbold_{\infty}$, and we will write
$\emptyset j = j$ for $j \in \Nbold$. 

In the applications mentioned above
the variables $\left(X_{\bi}\right)_{\bi \in \VV}$ 
of a RTP are often used as auxiliary variables to define or to construct some
useful random structures. 
In those cases typically the innovation process defines the
``internal'' variables while the RTP is constructed ``externally'' using
the consistency theorem. It is then natural
to ask whether the RTP is measurable only with respect to 
the i.i.d. innovation process $\left(\xi_{\bi}, N_{\bi}\right)$. 
\begin{Definition}
\label{Defi:Endogeny}
An invariant RTP with marginal $\mu$ is called \emph{endogenous}, if the
root variable $X_{\emptyset}$ is almost surely measurable with respect to the
$\sigma$-algebra 
\[ \GG := \sigma \left( \left\{ \left(\xi_{\bi}, N_{\bi} \right)
\,\Big\vert\, \bi \in \VV \,\right\} \right). \]
\end{Definition}
This notion of endogeny has been the main topic of discussion in
\cite{AlBa05}. The authors provide a necessary and sufficient
condition for endogeny in the general setup \cite[Theorem 11]{AlBa05}.
Some other concepts similar to endogeny can be found in \cite{Ba05b}. 

In this article we provide a non-trivial application of the theory
developed in \cite{AlBa05}. The example we consider here arise
from the study of the asymptotic limit of
\emph{random assignment problem} using local-weak convergence method 
\cite{Al01}. A detailed background of this example is given in 
Section \ref{Sec:Logistic-Intro}.

\subsection{Main Result}
\label{SubSec:Main-Result} 
The following RDE plays the central role in deriving the asymptotic
limit of the random assignment problem \cite{Al01}, 
\begin{equation}
X \ed \min_{j \geq 1} \left( \xi_j - X_j \right) \,\,\,\,
\mbox{on} \,\,\, \Rbold, 
\label{Equ:Logistic-RDE}
\end{equation}
where $\left(X_j\right)_{j \geq 1}$ are i.i.d with same law as $X$ and 
are independent of $\left(\xi_j\right)_{j \geq 1}$ which are
points of a Poisson point process of rate $1$ on $(0, \infty)$. 
It is known \cite{Al01} that the RDE (\ref{Equ:Logistic-RDE})
has a unique solution as the \emph{Logistic} distribution, given by
\begin{equation}
\bP\left( X \leq x \right) = \frac{1}{1 + e^{-x}}, \,\,\, x \in \Rbold.
\label{Equ:Logistic-defi} 
\end{equation} 
For this reason we will call RDE (\ref{Equ:Logistic-RDE}) the
\emph{Logistic RDE}. 
The following is our main result. 
\begin{Theorem}
\label{Thm:Logistic}
The invariant \emph{recursive tree process} with Logistic marginals
associated with the RDE (\ref{Equ:Logistic-RDE}) is endogenous. 
\end{Theorem}

This result though looks technical but, provides a concrete example falling
under the general theory developed in \cite{AlBa05}. 
The proof of Theorem \ref{Thm:Logistic}
involves analytic techniques, thus this work also demonstrate the need
of developing analytic tools for studying \emph{max-type} RDEs in general. 

\subsection{Outline of Rest of the Paper}
\label{SubSec:Outline}
The next section provides the background and motivation for deriving our
main result.   
In Section \ref{Sec:Review-Survey} we review some of the concepts from
\cite{AlBa05} and state a version of Theorem 11 of \cite{AlBa05},
which we will need to prove our main result. In 
Sections \ref{Sec:Logistic-biuni} and \ref{Sec:Logistic-Proof} we prove the
main result.
Finally Section \ref{Sec:Remarks} provides some further 
discussion. Some known facts about Logistic distribution which are needed
for the proofs are given in the appendix.

\section{Background and Motivation for Logistic RDE}
\label{Sec:Logistic-Intro} 
For a given $n \times n$ matrix of costs 
$\left(  C_{i j}  \right)$, 
consider the problem of assigning $n$ jobs to $n$ machines in the most 
``cost effective" way. Thus the task is to find a permutation $\pi$
of $\{ 1,2,\ldots,n \}$, which solves the following minimization
problem
\begin{equation} 
A_n := \min_{\pi} \sum_{i=1}^n C_{i, \pi(i)}.
\end{equation}
This problem has been extensively studied in literature for a fixed
cost matrix, and there are various algorithms to find the optimal permutation
$\pi$. A probabilistic model for the assignment problem
can be obtained by assuming that the costs are independent random
variables each with  Uniform$[0,1]$ distribution. Although 
this model appears to be quite simple, careful investigations of
it in the last few decades have shown that it has enormous
richness in its structure. See \cite{St97, AlSt04} 
for survey and other related works. 

Our interest in this problem is from another perspective. In 2001 Aldous
\cite{Al01} showed 
\begin{equation}
\lim_{n \rightarrow \infty} \bE[ A_n ] = \zeta(2) = \frac{\pi^2}{6},
\label{zeta-2-limit}
\end{equation}
confirming the earlier work of M\'{e}zard and Parisi \cite{MePa87}, where
they computed the same limit using some non-rigorous arguments based on
the \emph{replica method} \cite{MePaVi87}.
In an earlier work Aldous \cite{Al92a} showed that 
the limit of $\bE \left[ A_n \right]$ as
$n \rightarrow \infty$ exists for any i.i.d. cost distribution. He
also proved that the final limit does not
depend on the specifics of the cost distribution, 
except only on the value of the density at $0$, provided it exists and is 
strictly positive. So for
calculation of the limiting constant
one can assume that $C_{ij}$'s are i.i.d. with
Exponential distribution with mean $n$. Then we can redefine
the objective function $A_n$ in the normalized form, 
\begin{equation}
A_n := \min_{\pi} \frac{1}{n} \sum_{i=1}^n C_{i, \pi(i)}.
\label{Equ:def-objective-function}
\end{equation}
From historical perspective it is worth mentioning that in 1998 
Parisi \cite{Pa98} conjectured that in this case
the following exact formula holds
\[
\bE\left[A_n\right] = 1 + \frac{1}{4} + \cdots + \frac{1}{n^2},
\,\,\, \forall \,\,\, n \geq 1.
\]
Recently two separate groups
Linusson and W\"{a}stlund \cite{LW03} and 
C. Nair, B. Prabhakar and
M. Sharma \cite{NPS03} 
have independently proved this conjecture using combinatorial techniques.
Thus also proving the limit. 
However Aldous \cite{Al01} used local-weak convergence techniques to 
identify the limit constant
$\zeta(2)$ in terms of an optimal matching problem on an infinite
tree with random edge weights, described as follows
\begin{quote}
Let $\Tbold_{\infty} := \left( \VV, \EE \right)$ 
be the canonical infinite rooted labeled tree, as before, where
$\emptyset$ is the root. For every vertex 
$\bi \in \VV$, let $\left( \xi_{\bi j} \right)_{j \geq 1}$ be points
of a Poisson point process of rate $1$ on $(0,\infty)$, and they are
independent as $\bi$ varies. 
Define the weight of the edge $e=(\bi, \bi j) \in \EE$ as $\xi_{\bi j}$. 
\end{quote}
This structure is called \emph{Poisson weighted infinite tree} and henceforth
abbreviated as PWIT. 

Let $K_{n,n}^r$ be the complete graph on $n$ vertices with
a root selected uniformly at random. Suppose we also equip it with
i.i.d. Exponential edge weights with mean $n$. Then one can
show \cite{Al01, AlSt04} that in the sense of Aldous-Steel local weak
convergence $K_{n,n}^r$ converges to the PWIT. Moreover heuristically
the random assignment problem on $K_{n,n}^r$ has a ``natural'' analog to the
limit structure, which is to consider the ``optimal'' (in sense of
minimizing the ``total cost'') matching problem on PWIT. Naturally PWIT
being an infinite graph with edge weights each having mean at least $1$, 
the ``total cost'' of any matching is infinite a.s., and hence minimizing 
``total cost'' is not quite meaningful. However Aldous \cite{Al01}
showed that it is possible to make a sensible definition of
``optimal matching'' on PWIT which is invariant with respect to the
automorphism of the tree $\Tbold_{\infty}$, and minimizes the
``average edge weight''. This construction is quite hard, and we refer 
the readers to \cite{Al01, AlSt04} for the technical details. 
Here we only provide the basic essentials to understand the motivation
for our work. 

Consider the heuristic description of the ``optimal'' matching problem
on PWIT and suppose we define variables $X_{\bi}$ for each vertex
$\bi$ as follows 
\begin{eqnarray} 
X_{\bi} & = & \, \mbox{Total cost of a maximal matching on the subtree\ }
              \Tbold_{\infty}^{\bi} \nonumber \\
        &    & \!\!\!\!\! 
              - \, \mbox{Total cost of a maximal matching on the forest\ }
              \Tbold_{\infty}^{\bi} \setminus \{\bi\},
              \label{Equ:RTP-Heuristic}
\end{eqnarray} 
where $\Tbold_{\infty}^{\bi}$ is the subtree rooted at the vertex $\bi$. 
Here by ``total cost'' we mean the sum total of all the edge weights in 
the matching. As noted above, both the ``total costs'' appearing 
in (\ref{Equ:RTP-Heuristic}) are infinity almost surely. 
Thus rigorously speaking 
$X_{\bi}$ is not well defined. But at the heuristic level if we forget 
this important issue, and work with these $X_{\bi}$-variables as if they
are well defined, then simple manipulation yields that they must satisfy
the following recurrence relation (see Section 4.2 of \cite{Al01})
\begin{equation}
X_{\bi} = \min_{j \geq 1} \left( \xi_{\bi j} - X_{\bi j} \right) .
\label{Equ:Heuristic-RDE}
\end{equation}
This is of course the recurrence relation for a RTP associated with 
the Logistic RDE (\ref{Equ:Logistic-RDE}). 
Having observe that one can now construct the $X_{\bi}$-variables 
\emph{externally} as the RTP associated with the Logistic RDE, and
use them to redefine the optimal matching on PWIT. This is preciously
what Aldous did in \cite{Al01}, and later on referred as 
\emph{540-degree argument} by Aldous and Bandyopadhyay in \cite{AlBa05}. 
This construction also provides a 
characterization of the optimal matching on the PWIT. Finally one can then
derive the $\zeta(2)$-limit for the random assignment problem. 

Once again a natural question would be 
to figure out whether the random variables $X_{\bi}$'s are
truly external or not, in other words to see whether the RTP is endogenous 
or not (see remarks (4.2.d) and (4.2.e) in \cite{Al01}). This is our
main motivation for this work. Theorem \ref{Thm:Logistic} proves that
the $X_{\bi}$-variables can be defined using only the edge-weights and
hence they have no external randomness in them. 

Other significance
of this result has been pointed out in Section 7.5 of \cite{AlBa05}. 
We would like to note that
the endogeny of the Logistic RTP helps to define approximately feasible 
solution for the finite $n$-matching problem by 
using the optimal solution of the matching problem on PWIT. Thus 
with the help of endogeny one can write a possibly simpler
proof of Aldous' original argument for the $\zeta(2)$-limit of the
random assignment problem. But such derivation for this particular problem 
is not quite illuminating, and hence we do not pursue in that
direction. As indicated in Section 7.5 of \cite{AlBa05} in general
endogeny is an essential ingredient to make rigorous argument for 
the \emph{cavity method}, and this work is only to illustrate one such
non-trivial proof of endogeny.

\section{Review of Bivariate Uniqueness and Endogeny}
\label{Sec:Review-Survey}
In this section we review some of the concepts from \cite{AlBa05} which
will be needed to prove our main result, Theorem \ref{Thm:Logistic}. 

In the general setting of equation (\ref{Equ:RDE}) the question of 
endogeny is quite abstract. Aldous and Bandyopadhyay in \cite{AlBa05} 
introduces a concept called \emph{bivariate uniqueness} for an
invariant RTP, and showed under certain conditions that is equivalent to
endogeny. In the general setting bivariate uniqueness is defined as 
follows.

Consider a general RDE given by (\ref{Equ:RDE}) and let 
$T \colon \PP \rightarrow \PP\left(S\right)$ be the 
induced operator. We will consider a bivariate version of it. Write
$\PP^{(2)}$ for the space of probability measures
on $S^2 = S \times S$, with marginals in $\PP$. 
We can now define 
a map $T^{(2)} : \PP^{(2)} \rightarrow \PP\left(S^2\right)$ as follows
\begin{Definition}
\label{Defi:Bi-T-2nd}
For a probability $\mu^{(2)} \in \PP^{(2)}$, 
$T^{(2)}\left(\mu^{(2)}\right)$ is the joint distribution of 
\[
\left( \begin{array}{c}
       g\left( \xi, X_j^{(1)}, 1 \leq j \lleq N \right) \\
       g\left( \xi, X_j^{(2)}, 1 \leq j \lleq N \right)
       \end{array} \right)
\]
where we assume 
\begin{enumerate}
\item $\left(X_j^{(1)}, X_j^{(2)}\right)_{j \geq 1}$ are
      independent with joint distribution $\mu^{(2)}$ on $S^2$; 
\item the family of random variables 
      $\left(X_j^{(1)}, X_j^{(2)}\right)_{j \geq 1}$ are independent of
      the innovation pair $\left(\xi, N\right)$. 
\end{enumerate}
\end{Definition}
We note that we use the \emph{same realization} of the pair
$\left(\xi, N\right)$ in both components. Immediately from the definition
we have
\begin{itemize}
\item[(a)] If $\mu$ is a solution of the RDE then the associated
           \emph{diagonal measure} $\mudiag$ is a fixed-point for
           the operator $T^{(2)}$, where
           \begin{equation}
           \mudiag := \mbox{dist}\left(X, X\right) ,
           \label{Equ:Defi-mu-diag}
           \end{equation}
           where $X \sim \mu$. 
\item[(b)] If $\mu^{(2)}$ is a fixed-point of the operator $T^{(2)}$ then
           each marginal is a solution of the original RDE. 
\end{itemize}
So if $\mu$ is a solution of the RDE (\ref{Equ:RDE}) then $\mudiag$ is a
fixed point of $T^{(2)}$ and there may or may not be other fixed points
of $T^{(2)}$ with marginals $\mu$. 
\begin{Definition}
\label{Defi:Bi-Uni}
An invariant RTP with marginal $\mu$ has the \emph{bivariate uniqueness}
property if $\mudiag$ is the unique fixed point of the operator
$T^{(2)}$ with marginals $\mu$. 
\end{Definition}
Sometimes with slight abuse of terminology we will say that a solution 
of the RDE (\ref{Equ:RDE}) has bivariate uniqueness property, or even
the RDE has bivariate uniqueness property, if it has unique solution,
meaning that the invariant RTP associated with the solution has the
bivariate uniqueness property. Similar abuse will be done for the
term endogeny also. 

Theorem 11 of \cite{AlBa05} shows that under appropriate assumptions the
two concepts, namely bivariate uniqueness and endogeny are equivalent. 
Rather than stating this general equivalence theorem, we here only state
the part we will need to prove endogeny for the Logistic RDE. 

\begin{Theorem}[Theorem 11(b) of \cite{AlBa05}]
\label{Thm:Equiv-Part-b}
Let $S$ be a Polish space. Consider an invariant RTP with marginal
distribution $\mu$. Suppose the bivariate uniqueness property holds. 
If also $T^{(2)}$ is continuous with respect to the weak convergence on 
the set of bivariate distributions with marginals $\mu$, then the
endogenous property holds. 
\end{Theorem} 

Thus to prove endogenous property for the Logistic RDE (\ref{Equ:Logistic-RDE})
we will show that the bivariate uniqueness property holds and also
establish the technical condition of Theorem \ref{Thm:Equiv-Part-b}, these
are done in the following two sections.

\section{Bivariate Uniqueness for the Logistic RDE}
\label{Sec:Logistic-biuni} 
In this section we prove the bivariate uniqueness property for the
Logistic RDE (\ref{Equ:Logistic-RDE}). 

\begin{Theorem}
\label{Thm:bivariate}
Consider the following bivariate RDE
\begin{equation}
\left( \begin{array}{c} X \\ Y \end{array} \right)
\ed 
\left( \begin{array}{c} 
       \mathop{\min}\limits_{j \geq 1} \left( \xi_j - X_j \right) \\
       \mathop{\min}\limits_{j \geq 1} \left( \xi_j - Y_j \right)
       \end{array} \right),
\label{bivariate-eqn}
\end{equation}
where $\left( X_j, Y_j \right)_{j \geq 1}$ are i.i.d. pairs with
same joint distribution as $\left( X, Y \right)$ and are independent of
$\left( \xi_j \right)_{j \geq 1}$ which are points of a Poisson 
process of rate $1$ on $(0,\infty)$. 
Then the unique solution of this RDE is given by 
the diagonal measure 
$\mudiag$ where $\mu$ is the Logistic distribution.
\end{Theorem}  

\subsection{Proof of Theorem \ref{Thm:bivariate}} 
\label{Subsec:Logistic-biuni-Proof}

First observe that if the equation (\ref{bivariate-eqn}) has a
solution then, the marginal distributions of $X$ and $Y$ solve
the Logistic RDE (\ref{Equ:Logistic-RDE}), 
and hence they are both Logistic. Further by inspection 
$\mudiag$ is a solution of (\ref{bivariate-eqn}). So it is enough to
prove that $\mudiag$ is the only solution of (\ref{bivariate-eqn}). 

Let $\mu^{(2)}$ be a solution of (\ref{bivariate-eqn}). Notice that the points 
$\left\{ \left( \xi_j ; (X_j, Y_j) \right) \left\vert\right. j \geq 1 \right\}$
form a Poisson point process, say ${\mathscr P}$, 
on $(0,\infty) \times \Rbold^2$,
with mean intensity 
$\rho(t;(x,y)) \, dt \, d(x,y) := dt \, \mu^{(2)}(d(x,y))$.
Thus if $G(x,y) := \bP \left( X > x, Y>y \right)$, for 
$x, y \in \Rbold$, then 
\begin{eqnarray}
G(x,y) & = & \bP \left( \min_{j \geq 1} \left( \xi_j - X_j \right) > x, 
             \,\, \mbox{and,} \,\,\min_{j \geq 1} \left( \xi_j - Y_j
             \right) > y \right) \nonumber \\
 & = & \bP \left( \mbox{No points of\ } {\mathscr P} \mbox{\ are in\ } 
       \left\{ (t;(u,v)) \, \Big\vert \, t-u \leq x, \mbox{\ or,\ } 
       t-v \leq y \right\} \right) \nonumber \\
 & = & \exp \left( - 
       \mathop{\int\int\int}\limits_{t-u \leq x,
             \mbox{\ or,\ } t-v \leq y} \! \rho(t;(u,v)) \, dt \, d(u,v)
       \right) \nonumber \\
 & = & \exp\left(- \int_0^{\infty} \! 
       \left[ \Hbar(t-x) + \Hbar(t-y) - G(t-x, t-y) \right] \, dt
       \right) \nonumber \\
 & = & \Hbar(x) \, \Hbar(y) \, \exp\left( \int_0^{\infty} \!
       G(t-x, t-y) \, dt \right), \label{G-eqn}
\end{eqnarray}
where $\Hbar$ is the right tail of Logistic distribution, defined as
$\Hbar(x) = e^{-x}/\left(1+e^{-x}\right)$ for $x \in \Rbold$.   
The last equality follows from properties of the Logistic distribution 
(see Fact \ref{fact:H-3} of appendix). 
For notational convenience in this paper we will write
${\overline F}\left( \cdot \right) := 1 - F\left( \cdot \right)$, for any
distribution function $F$.

The following simple lemma reduces the bivariate problem to a
univariate problem. 
\begin{Lemma}
\label{lem:min}
For any two random variables $U$ and $V$, $U=V$ a.s. if and only
if $U \ed V \ed U \wedge V$. 
\end{Lemma}
  
\noindent
\emph{Proof :} First of all if $U=V$ a.s. then $U \wedge V = U$ a.s. 

Conversely suppose that $U \ed V \ed U \wedge V$. Fix a rational 
$q$, then under our assumption, 
\begin{eqnarray*}
\bP \left( U \leq q < V \right) & = & \bP \left( V > q \right) -
                                      \bP \left( U > q, \, V > q
                                      \right) \\
                                & = & \bP \left( V > q \right) -
                                      \bP \left( U \wedge V > q
                                      \right) \\
                                & = & 0
\end{eqnarray*}
A similar calculation will show that 
$\bP \left( V \leq q < U \right) = 0$. These are true for
any rational $q$, thus $\bP \left( U \neq V \right) = 0$. 
\qed  

Thus if we can show that $X \wedge Y$ also has Logistic distribution,
then from the lemma above we will be able to conclude that $X = Y$
a.s., and hence the proof will be complete. Put 
$g(\cdot) := \bP \left( X \wedge Y > \cdot \right)$, we will show
$ g = \Hbar$.  Now, for every fixed $x \in \Rbold$, by definition 
$g(x) = G(x,x)$. So using (\ref{G-eqn}) we
get
\begin{equation} 
g(x) = \Hbar^2(x) \, \exp\left( \int_{-x}^{\infty} \! g(s) \, ds \right),
\,\, x \in \Rbold.
\label{g-eqn}
\end{equation}
Notice that from (A1) (see Fact \ref{fact:H-3} of appendix) 
$g = \Hbar$ is a solution of this non-linear 
integral equation (\ref{g-eqn}), which corresponds to the solution
$\mu^{(2)} = \mudiag$ of the original equation (\ref{bivariate-eqn}). 
To complete the proof of Theorem
\ref{Thm:bivariate} we need to show that this is
the only solution. For that we will prove
that the operator associated with (\ref{g-eqn}) (defined on an
appropriate space) is monotone and has unique fixed-point as
$\Hbar$. The techniques we will use here are similar to Eulerian
recursion \cite{Sim72}, and are heavily based on analytic arguments. 

Let ${\mathfrak F}$ be the set of all functions
$f : \Rbold \rightarrow [0,1]$ such that 
\begin{itemize}
\item $\Hbar^2(x) \leq f(x) \leq \Hbar(x), \,\, \forall \,\, x \in \Rbold$, 
\item $f$ is continuous and non-increasing. 
\end{itemize}
Observe that by definition $\Hbar \in {\mathfrak F}$. Further from 
(\ref{g-eqn}) it follows that 
$g(x) \geq \Hbar^2(x)$, as well as, 
$g(x) = \bP \left( X \wedge Y > x \right) \leq \bP \left( X > x \right)
= \Hbar(x), \,\, \forall \,\, x \in \Rbold$. Note also that 
$g$ being the tail of the random variable $X \wedge Y$, is continuous 
(because both $X$ and $Y$ are continuous random variables) and 
non-increasing. So it is appropriate to
search for solutions of (\ref{g-eqn}) in ${\mathfrak F}$. 

Let $T : {\mathfrak F} \rightarrow {\mathfrak F}$ be defined as 
\begin{equation} 
T(f)(x) := \Hbar^2(x) \, \exp\left( \int_{-x}^{\infty} \! f(s) \, ds
\right), \,\, x \in \Rbold.
\label{def-T}
\end{equation}
Note that this operator $T$ is not same as the general operator defined in
Section \ref{Sec:Intro}, henceforth by $T$ we will mean the specific
operator defined above. 
Proposition \ref{tech-1} of Section \ref{Subsec:Technical} shows that $T$
does indeed map ${\mathfrak F}$ into itself. 
Observe that the equation (\ref{g-eqn}) is
nothing but the fixed-point equation associated with the operator $T$,
that is, 
\begin{equation}
g = T(g) \,\,\, \mbox{on\ \ } {\mathfrak F}.
\label{T-eqn}
\end{equation}
We here note that using (A1) (see Fact \ref{fact:H-3} of appendix) $T$ 
can also be written as
\begin{equation}
T(f)(x) := \Hbar(x) \, \exp\left( - \int_{-x}^{\infty} \! 
\left( \Hbar(s)- f(s) \right) \, ds \right), \,\, x \in \Rbold,
\label{working-def-T}
\end{equation} 
which will be used in the subsequent discussion. 

Define a partial order $\stleq$ on 
${\mathfrak F}$ as, $f_1 \stleq f_2$ in
${\mathfrak F}$ if $f_1(x) \leq f_2(x), \,\, \forall \,\, x \in \Rbold$, 
then the following result holds.
\begin{Lemma} 
\label{lem:T-monotone}
$T$ is a monotone operator on the partially ordered set 
$({\mathfrak F}, \stleq)$.
\end{Lemma}

\noindent
\emph{Proof :} Let $f_1 \stleq f_2$ be two elements of ${\mathfrak F}$, 
so from definition
$f_1(x) \leq f_2(x), \,\, \forall \,\, x \in \Rbold$. Hence
\[
\begin{array}{ccccc}
  & \mathop{\int}\limits_{-x}^{\infty} \! f_1(s) \, ds & \leq & 
\mathop{\int}\limits_{-x}^{\infty} \! f_2(s) \, ds, & 
\forall \,\, x \in \Rbold \\
\Rightarrow & T(f_1)(x) & \leq & T(f_2)(x), & \forall \,\, x \in \Rbold \\
\Rightarrow & T(f_1) & \stleq & T(f_2). &   
\end{array}
\]
\qed

Put $f_0 = \Hbar^2$, and for $n \in \Nbold$, define $f_n \in {\mathfrak F}$
recursively as, $f_n = T(f_{n-1})$. Now from Lemma \ref{lem:T-monotone} we get
that if $g$ is a fixed-point of $T$ in ${\mathfrak F}$ then, 
\begin{equation}
f_n \stleq g, \,\,\, \forall \,\,\, n \geq 0. 
\label{lower-bound}
\end{equation}
If we can show $f_n \rightarrow \Hbar$ pointwise, then using 
(\ref{lower-bound}) we will get $\Hbar \stleq g$, so from
definition of ${\mathfrak F}$ it will follow that  
$g = \Hbar$, and our proof will be complete. 
For that, the following lemma gives an explicit recursion for the
functions $\left\{ f_n \right\}_{n \geq 0}$. 
\begin{Lemma}
\label{lem:recurssion} 
Let $\beta_0(s) = 1 - s$, $0 \leq s \leq 1$. Define recursively
\begin{equation}
\beta_n(s) := \int_s^1 \! \frac{1}{w} \left( 1 - e^{-\beta_{n-1}(1-w)}
\right) \, dw, \,\, 0 < s \leq 1. 
\label{def-beta}
\end{equation}
Then for $n \geq 1$, 
\begin{equation}
f_n (x) = \Hbar(x) \, \exp\left( - \beta_{n-1}(\Hbar(x)) \right),
\,\, x \in \Rbold.
\label{f-recurssion}
\end{equation}
\end{Lemma}

\noindent
\emph{Proof :} We will prove this by induction on $n$. Fix $x \in \Rbold$, for 
$n=1$ we get 
\begin{eqnarray*}
f_1 (x) & = & T(f_0)(x) \\
        & = & \Hbar(x) \, \exp\left( - \int_{-x}^{\infty} \!
              \left( \Hbar(s) - \Hbar^2(s) \right) \, ds \right) 
              \mbox{\ \ \ \ \ [using (\ref{working-def-T})]} \\
        & = & \Hbar(x) \, \exp\left( - \int_{-x}^{\infty} \!
              \Hbar(s) \, \left( 1 - \Hbar(s) \right) \, ds \right) \\
        & = & \Hbar(x) \, \exp\left( - \int_{-x}^{\infty} \!
              \Hbar(s) \, H(s) \, ds \right) \\
        & = & \Hbar(x) \, \exp\left( - \int_{-x}^{\infty} \!
              H'(s) \, ds \right) 
              \mbox{\ \ \ \ \ [using Fact \ref{fact:H-1} of appendix]} \\
        & = & \Hbar(x) \, \exp\left( - H(x) \right) \\
        & = & \Hbar(x) \, \exp\left( - \beta_0(\Hbar(x)) \right)
\end{eqnarray*}

Now, assume that the assertion of the Lemma is true for
$n \in \{1,2,\ldots, k\}$, for some $k \geq 1$, then from definition
we have 
\begin{eqnarray}
f_{k+1}(x) & = & T(f_k)(x) \nonumber \\ 
           & = & \Hbar(x) \, \exp\left( - \int_{-x}^{\infty} \!
                 \left( \Hbar(s) - f_k (s) \right) \, ds \right)
                 \mbox{\ \ \ \ \ [using (\ref{working-def-T})]}
                 \nonumber \\
           & = & \Hbar(x) \, \exp\left( - \int_{-x}^{\infty} \!
                 \Hbar(s) \left( 1 - e^{-\beta_{k-1}(\Hbar(s))} \right)
                 \, ds \right) \nonumber \\
           & = & \Hbar(x) \, \exp\left( - \int_{\Hbar(x)}^1 \!
                 \frac{1}{w} \left( 1 - e^{-\beta_{k-1}(1-w)} \right)
                 \, dw \right) \label{temp1}
\end{eqnarray}
The last equality follows by substituting $w=H(s)$ and thus from 
Fact \ref{fact:H-1} and Fact \ref{fact:H-2} of the appendix
we get that $\sfrac{dw}{w} = \Hbar(s) \, ds$ and $H(-x) = \Hbar(x)$. 
Finally by definition of $\beta_n$'s and using (\ref{temp1}) we get 
$f_{k+1} = T(f_k)$. \qed 
  
To complete the proof it is now enough to show that $\beta_n \rightarrow 0$
pointwise, which will imply by Lemma \ref{lem:recurssion} that 
$f_n \rightarrow \Hbar$ pointwise, as $n \rightarrow \infty$. 
Using Proposition \ref{tech-2} (see Section \ref{Subsec:Technical})
we get the following 
characterization of the pointwise limit of these $\beta_n$'s. 
\begin{Lemma} 
\label{lem:L}
There exists a function $L : [0,1] \rightarrow [0,1]$ with $L(1)=0$,
such that 
\begin{equation} 
L(s) = \int_s^1 \! \frac{1}{w} \left( 1 - e^{-L(1-w)} \right) \, dw, 
\,\, \forall \, s \in [0,1), 
\label{L-eqn}
\end{equation}
and $L(s) = \mathop{\lim}\limits_{n \rightarrow \infty} \beta_n(s), 
\,\, \forall \,\, 0 \leq s \leq 1$. 
\end{Lemma}

\noindent
\emph{Proof :} From the Proposition \ref{tech-2} we know that for any 
$s \in [0,1]$ the sequence $\left\{ \beta_n(s) \right\}$ is
decreasing, and hence
$\exists$ a function $L : [0,1] \rightarrow [0,1]$ such that 
$L(s) = \mathop{\lim}\limits_{n \rightarrow \infty} \beta_n(s)$. 
Now observe that 
$\beta_n(1-w) \leq \beta_0(1-w) = w, \,\, \forall \,\, 0 \leq w \leq 1$,
and hence 
\[
0 \leq \frac{1}{w} \left( 1 - e^{-\beta_n(1-w)} \right) \leq 
\frac{\beta_n(1-w)}{w} \leq 1, \,\,\, \forall \,\,\, 0 \leq w \leq 1.
\]
Thus by taking limit as $n \rightarrow \infty$ in (\ref{def-beta}) and
using the \emph{dominated convergence theorem}
along with part (a) of Proposition \ref{tech-2} we get that
\[
L(s) = \int_s^1 \! \frac{1}{w} \left( 1 - e^{-L(1-w)} \right)
      \, dw, \,\,\, \forall \,\,\, 0 \leq s < 1.
\]
\qed

The above lemma basically translates the non-linear integral equation
(\ref{g-eqn}) to the non-linear integral equation (\ref{L-eqn}), where the
solution $g=\Hbar$ of (\ref{g-eqn}) is given by the solution 
$L \equiv 0$ of (\ref{L-eqn}). So at first sight
this may not lead us to the conclusion. But fortunately, something
nice happens for equation (\ref{L-eqn}), and we have the following result
which is enough to complete the proof of Theorem
\ref{Thm:bivariate}.

\begin{Lemma} 
\label{lem:L-problem}
If $L : [0,1] \rightarrow [0,1]$ is a function which
satisfies the non-linear integral equation (\ref{L-eqn}), namely, 
\[
L(s) = \int_s^1 \! \frac{1}{w} \left( 1 - e^{-L(1-w)} \right) \, dw,
\,\,\, \forall \,\,\, 0 \leq s < 1,
\]
and if $L(1) = 0$, then $L \equiv 0$.
\end{Lemma}

\noindent
\emph{Proof :} First note that $L \equiv 0$ is a solution. 
Now let $L$ be any solution of (\ref{L-eqn}), then $L$ is infinitely
differentiable on the open interval $(0,1)$, by repetitive application
of Fundamental Theorem of Calculus. 

Consider, 
\begin{equation}
\eta(w) := (1-w) e^{L(1-w)} + w e^{-L(w)} - 1, \,\, w \in [0,1]. 
\label{def-eta}
\end{equation}
Observe that $\eta(0)=\eta(1)=0$ as $L(1)=0$. Now, from (\ref{L-eqn})
we get that 
\begin{equation}
L'(w) = - \frac{1}{w} \left( 1 - e^{-L(1-w)} \right), \,\, w \in (0,1).
\label{L-derivative}
\end{equation}
Thus differentiating the function $\eta$ we get 
\begin{equation}
\eta'(w) = e^{-L(w)} \left[ 2 - \left( e^{L(1-w)} + e^{-L(1-w)} \right)
\right] \leq 0, \,\, \forall \,\, w \in (0,1).
\label{eta-derivative}
\end{equation}
So the function 
$\eta$ is decreasing in $(0,1)$ and is continuous in $[0,1]$ with
boundary values as $0$, hence $\eta \equiv 0$ Thus we must have
$\eta' \equiv 0$, so from equation (\ref{eta-derivative}) we get that
\[
e^{L(s)} + e^{-L(s)} = 2 \,\,\,\, \mbox{for all} \,\,\, s \in (0,1).
\]
This implies $L \equiv 0$ on $[0,1]$. \qed

\subsection{Some Technical Details}
\label{Subsec:Technical}

This section provides some of the technical results which were needed 
in the previous section. 

\begin{Proposition}
\label{tech-1}
The operator $T$ maps ${\mathfrak F}$ into ${\mathfrak F}$.
\end{Proposition} 

\noindent
\emph{Proof :} First note that if 
$f \in {\mathfrak F}$, then by definition 
$T(f)(x) \geq \Hbar^2(x), \,\, \forall \,\, x \in \Rbold$. 
Next by definition of ${\mathfrak F}$ we get that 
$f \in {\mathfrak F} \Rightarrow f \stleq \Hbar$, thus
\begin{eqnarray*}
          &             & \int_{-x}^{\infty} \! f(s) \, ds \,\,
                          \leq \, \int_{-x}^{\infty} \! \Hbar(s) \, ds,
                          \,\, \forall \,\, x \in \Rbold \\
          & \Rightarrow & T(f)(x) \,\, \leq \, \Hbar^2(x) \, 
                          \exp\left( \int_{-x}^{\infty} \! \Hbar(s) \, ds
                          \right) \,\, = \,\, \Hbar(x), 
                          \,\, \forall \,\, x \in \Rbold 
\end{eqnarray*}
The last equality follows from (A1) (see Fact \ref{fact:H-3} of
appendix). So, 
\begin{equation} 
\Hbar^2(x) \leq T(f)(x) \leq \Hbar(x), \,\,\, \forall \,\,\, x \in \Rbold.
\label{tech-11}
\end{equation}

Now we need to show that for any $f \in {\mathfrak F}$ we must have $T(f)$ 
continuous and non-increasing. 
From the definition $T(f)$ is continuous (in fact, infinitely
differentiable). Moreover if $x \leq y$ be two real numbers, then
\[
\int_{-x}^{\infty} \! \left( \Hbar(s) - f(s) \right) \, ds
\leq \int_{-y}^{\infty} \! \left( \Hbar(s) - f(s) \right) \, ds, 
\]
because $f \stleq \Hbar$. Also $\Hbar(x) \geq \Hbar(y)$, thus
using (\ref{working-def-T}) we get 
\begin{equation}
T(f)(x) \geq T(f)(y) 
\label{tech-13}
\end{equation}
So using (\ref{tech-11}) and (\ref{tech-13}) we conclude
that $T(f) \in {\mathfrak F}$ if $f \in {\mathfrak F}$. 
\qed

\begin{Proposition}
\label{tech-2}
The following are true for the sequence of functions 
$\left\{ \beta_n \right\}_{n \geq 0}$ defined in (\ref{def-beta}).
\begin{itemize}
\item[(a)] For every fixed $s \in (0,1]$, the sequence 
           $\left\{ \beta_n(s) \right\}$ is decreasing.
\item[(b)] For every $n \geq 1$, 
           $\mathop{\lim}\limits_{s \rightarrow 0+} \beta_n(s)$ exists,
           and is given by 
           \[ \int_0^1 \! \sfrac{1}{w} \left( 1 - e^{-\beta_{n-1}(1-w)} 
           \right) \, dw, \]
           we will write this as $\beta_n(0)$.
\item[(c)] The sequence of numbers
           $\left\{ \beta_n(0) \right\}$ is also decreasing. 
\end{itemize}
\end{Proposition}

\noindent
\emph{Proof :} (a) Notice that $\beta_0(s) = 1 - s$ for $s \in [0,1]$, thus
\[ 
\beta_1(s) = \int_s^1 \! \frac{1 - e^{-w}}{w} \, dw
< 1 - s = \beta_0(s), \,\,\, \forall \,\,\, s \in (0,1].
\]
Now assume that for some $n \geq 1$ we have 
$\beta_n(s) \leq \beta_{n-1}(s) \leq \cdots \leq \beta_0(s), \,\,\,
\forall \,\, s \in (0,1]$, if we show that 
$\beta_{n+1}(s) \leq \beta_n(s), \,\,\, \forall \,\,\, s \in (0,1]$ then by
induction the proof will be complete. For that, fix $s \in (0,1]$ then 
\begin{eqnarray*}
\beta_{n+1}(s) & = & \int_s^1 \! \frac{1}{w} \left( 1 - e^{-\beta_n(1-w)}
                     \right) \, dw \\
           & \leq  & \int_s^1 \! \frac{1}{w} \left( 1 - e^{-\beta_{n-1}(1-w)}
                     \right) \, dw \\
               & = & \beta_n(s). 
\end{eqnarray*}
This proves the part (a). 

\vspace{0.15in}

\noindent
(b, c) First note that by trivial induction
$\beta_n(s) \geq 0$ for every $s \in (0,1]$, $n \geq 0$. Thus from
definition for every $n \geq 0$, the limit
$\mathop{\lim}\limits_{s \rightarrow 0+} \beta_n(s)$ exists in 
$[0, \infty]$ and is given by 
\[ \int_0^1 \! \sfrac{1}{w} \left( 1 - e^{-\beta_{n-1}(1-w)} 
           \right) \, dw. 
\]
Now using (a) above we conclude
\begin{equation}
\beta_{n+1}(0)   =  \lim_{s \rightarrow 0+} \beta_{n+1}(s)
               \leq \lim_{s \rightarrow 0+} \beta_n(s)
                 =  \beta_n(0) ,
\label{Equ:Beta-0-decreasing}
\end{equation}
for every $n \geq 0$. Since $\beta_0(0) = 1$, so we get 
$\beta_n(0) < \infty$ for all $n \geq 0$, and the sequence is decreasing. 
Proving parts (b) and (c). \qed

\section{Proof of Theorem \ref{Thm:Logistic}}
\label{Sec:Logistic-Proof} 

Once again we will use the general Theorem 11(b) of \cite{AlBa05}, stated here
as Theorem \ref{Thm:Equiv-Part-b}. 
We note that by 
Theorem \ref{Thm:bivariate} the Logistic RDE (\ref{Equ:Logistic-RDE}) has
bivariate uniqueness property and hence all remains is to check the
technical continuity condition. 

\begin{Proposition}
\label{prop:weak-continuity}
Let ${\mathfrak S}$ be the set of all probabilities on $\Rbold^2$ and let 
$\Gamma : {\mathfrak S} \rightarrow {\mathfrak S}$ be the operator
associated with the RDE (\ref{bivariate-eqn}), that is, 
\begin{equation}
\Gamma \left( \mu^{(2)} \right) 
\ed 
\left( \begin{array}{c} 
              \mathop{\min}\limits_{j \geq 1} \left( \xi_j - X_j \right) \\
              \mathop{\min}\limits_{j \geq 1} \left( \xi_j - Y_j \right)
       \end{array} \right), 
\label{Equ:Logistic-biop}
\end{equation}
where $\left( X_j, Y_j \right)_{j \geq 1}$ are i.i.d with joint law
$\mu^{(2)} \in {\mathfrak S}$ and are independent of 
$\left( \xi_j \right)_{j \geq 1}$ which are points of a Poisson point process
of rate $1$ on $(0,\infty)$.  
Then $\Gamma$ is continuous with respect to the weak convergence
topology when restricted
to the subspace ${\mathfrak S}^*$ defined as 
\begin{equation}
{\mathfrak S}^* := 
\left\{ \mu^{(2)} \, \Big\vert \, \mbox{both the marginals of\ }
\mu^{(2)} \mbox{\ are Logistic distribution}\, \right\}.
\label{Equ:Logistic-subspace}
\end{equation}
\end{Proposition} 

Before we prove this proposition, 
it is worth mentioning that the operator $\Gamma$ is not continuous 
with respect to the weak convergence topology on the
whole space ${\mathfrak S}$. In fact, as it turns out 
it is every where discontinuous on ${\mathfrak S}$ 
(see Section \ref{Sec:Remarks}). 
But fortunately for applying Theorem \ref{Thm:Equiv-Part-b} we only need
the continuity of $\Gamma$ when restricted to the subspace ${\mathfrak S}^*$.

\noindent
\emph{Proof of Proposition \ref{prop:weak-continuity} :} 
Let $\left\{ \mu^{(2)}_n \right\}_{n=1}^{\infty} \subseteq 
{\mathfrak S}^{\star}$ and suppose that 
$\mu^{(2)}_n \cd \mu^{(2)} \in {\mathfrak S}^{\star}$. We will show that 
$\Gamma(\mu^{(2)}_n) \cd \Gamma(\mu^{(2)})$.

Let $\left( \Omega, {\mathcal F}, \bP \right)$ be a probability space
such that, $\exists$ $\left\{ \left( X_n, Y_n \right) \right\}_{n=1}^{\infty}$
and $\left( X, Y \right)$ random vectors taking values in $\Rbold^2$, with
$\left( X_n, Y_n \right) \sim \mu^{(2)}_n, \,\, n \geq 1$, and 
$\left( X, Y \right) \sim \mu^{(2)}$. Notice that by definition  
$X_n \ed Y_n \ed X \ed Y$, and each has Logistic distribution. 

Fix $x,y \in \Rbold$, then using similar calculations as in  
(\ref{G-eqn}) we get
\begin{eqnarray} 
G_n(x, y) 
& := & \Gamma(\mu^{(2)}_n) \left( (x,\infty) \times (y,\infty) \right)
\nonumber \\
& = & \Hbar(x) \Hbar(y) \, \exp\left( - \int_0^{\infty} \! 
\bP\left( X_n > t-x, Y_n > t-y \right) \, dt \right) \nonumber \\
& = & \Hbar(x) \Hbar(y) \, \exp\left( - \int_0^{\infty} \! 
\bP\left( (X_n + x) \wedge (Y_n + y) > t \right) \, dt \right)
\nonumber \\
 & = & \Hbar(x) \Hbar(y) \, \exp\left( - \bE 
\left[ (X_n + x)^+ \wedge (Y_n + y )^+ \right] \right), \label{def-Gn} 
\end{eqnarray}
and a similar calculation will also give that 
\begin{eqnarray} 
G(x,y) 
& := & \Gamma(\mu^{(2)}) \left( (x,\infty) \times (y,\infty) \right) 
\nonumber \\
& = & \Hbar(x) \Hbar(y) \, \exp\left( - \bE 
\left[ (X + x)^+ \wedge (Y + y )^+ \right] \right). \label{def-G}
\end{eqnarray} 
Now to complete the proof all we need is to show
\[
\bE \left[ (X_n + x)^+ \wedge (Y_n + y )^+ \right]
\longrightarrow
\bE \left[ (X + x)^+ \wedge (Y + y )^+ \right].
\]
Since we assumed that $\left(X_n, Y_n \right) \cd \left( X, Y \right)$
thus 
\begin{equation} 
\left( X_n + x \right)^+ \wedge \left( Y_n + y \right)^+
\cd
\left( X + x \right)^+ \wedge \left( Y + y \right)^+,
\,\,\, \forall \,\,\, x,y \in \Rbold.
\label{conv-of-min}
\end{equation}

Fix $x,y \in \Rbold$, define 
$Z_n^{x,y} := \left( X_n + x \right)^+ \wedge \left( Y_n + y \right)^+$, 
 and $Z^{x,y} := \left( X + x \right)^+ \wedge \left( Y + y \right)^+$.
Observe that
\begin{equation} 
0 \, \leq \, Z_n^{x,y} \, \leq \, \left( X_n + x \right)^+ \leq \,
\left\vert X_n + x \right\vert, \,\,\, \forall \,\,\, n \geq 1.
\label{bound-Zn}
\end{equation}
But, $\left\vert X_n + x \right\vert \ed \left\vert X + x \right\vert, \,\,
\forall \,\, n \geq 1$. So clearly 
$\left\{ Z_n^{x,y} \right\}_{n=1}^{\infty}$ is uniformly
integrable. Hence we conclude (using Theorem 25.12 of Billingsley
\cite{Bill95}) that 
\[
\bE \left[ Z_n^{x,y} \right] 
\longrightarrow 
\bE \left[ Z^{x,y} \right].
\]
This completes the proof. \qed

\section{Final Remarks}
\label{Sec:Remarks}

(a) Intuitively, a natural approach to show that the fixed-point
equation $\Gamma(\mu^{(2)}) = \mu^{(2)}$ on ${\mathfrak S}$ 
has unique solution,
would be to specify a metric $\rho$ on ${\mathfrak S}$ such that the operator
$\Gamma$ becomes a contraction with respect to it. Unfortunately,
this approach seems rather hard or may even be impossible. Perhaps the
reason being the Logistic RDE (\ref{Equ:Logistic-RDE}) itself does not
have a contractive property, in fact, it does not
have a full domain of attraction (see \cite{AlBa05}). However 
its exact domain of attraction is not yet known 
(see open problem 62 of \cite{AlBa05}). On the other hand
from the proof of Theorem \ref{Thm:bivariate} it is clear that
equation (\ref{g-eqn}) has the whole of ${\mathfrak F}$ within its domain
of attraction. So it is possible to have a suitable metric of contraction
for $T$ but, we have been unable to find it. 

\vspace{0.25in} 

\noindent
(b) Although at first glance it seems that the operator $T$ as defined
in (\ref{def-T}) is just an analytic tool to solve the equation
(\ref{g-eqn}) but,  it has a nice interpretation through Logistic RDE
(\ref{Equ:Logistic-RDE}). Suppose ${\mathfrak A}$ is the operator associated
with Logistic RDE, that is, 
\begin{equation}
{\mathfrak A}(\mu) \ed \min_{j \geq 1} \left( \xi_j - X_j \right),  
\label{def-Logistic-operator} 
\end{equation} 
where $\left( \xi_j \right)_{j \geq 1}$ are points of a Poisson point process
of mean intensity $1$ on $(0,\infty)$, and are independent of 
$\left( X_j \right)_{j \geq 1}$, which are i.i.d with distribution
$\mu$ on $\Rbold$. It is easy to check that the domain of definition
of ${\mathfrak A}$ is the space
\begin{equation}
{\mathcal A} := \left\{ F \, \Big\vert \, 
F \mbox{\ is a distribution function on\ } \Rbold \mbox{\ and\ } 
\int_0^{\infty} \! {\overline F} (s) \, ds < \infty \right\}.
\label{def-domian-of-A}
\end{equation}
Note that the condition 
$\int_0^{\infty} \! {\overline F}(s) \, ds < \infty$ means 
$\bE_F \left[ X^+ \right] < \infty$. 
Now it is easy to see that ${\mathfrak F}$ can be embedded into 
${\mathcal A}$ and definition of $T$ can be naturally extended 
on whole of ${\mathcal A}$. In that case the following identity holds 
\begin{equation} 
\frac{\overline{T(\mu)}(\cdot)}{\Hbar(\cdot)} \times
\frac{\overline{{\mathfrak A}(\mu)}(\cdot)}{\Hbar(\cdot)} = 1, 
\,\,\, \forall \,\,\, \mu \in {\mathcal A}.
\label{identity-eqn}
\end{equation} 
This at least explains the monotonicity of $T$ through
anti-monotonicity property of the Logistic operator
${\mathfrak A}$ (easy to check).

\vspace{0.25in}

\noindent
(c) It is interesting to note that the operator ${\mathfrak A}$ 
is every where discontinuous with respect to the weak 
convergence topology on ${\mathcal A}$. 
This is because, given any distribution 
$F_0 \in {\mathcal A}$, we can construct a sequence
of distributions $\left\{ F_n \right\}_{n \geq 1} \subseteq {\mathcal A}$
converging in distribution to $F_0$, such that  
\[
\bE_{F_n}\left[ \left(x+X_n\right)^+ \right] \rightarrow \infty, 
\,\,\,\, \mbox{for all} \,\,\, x \in \Rbold . 
\]
Note $F_0 \in {\mathcal A} \, \Rightarrow \, 
\bE_{F_0}\left[ \left(x+X_0\right)^+ \right] < \infty$, for all 
$x \in \Rbold$. 
On the other hand we know that for any distribution function $F$, 
\[
{\mathfrak A}\left( F \right) \left( x \right) = 
1 - \exp \left( - \bE_{F}\left[ \left(x+X\right)^+ \right] \right), 
\,\,\,\, \mbox{for all} \,\,\, x \in \Rbold  
\]
(see the proof of Fact \ref{fact:H-3} in the appendix). 
Thus for every $x \in \Rbold$, 
\[ 
{\mathfrak A}\left(F_n\right)(x) \nrightarrow 
{\mathfrak A}\left(F_0\right)(x). 
\]
So ${\mathfrak A}$ is discontinuous at $F_0$ for every $F_0 \in {\mathcal A}$
with respect to the weak convergence topology. 
This also indicates that the same phenomenon is true for the
bivariate operator $\Gamma$.

\section*{Appendix}
\label{sec:appendix}

Here we provide some known facts about the Logistic distribution
which are used in the Sections \ref{Sec:Logistic-biuni} and 
\ref{Sec:Logistic-Proof}. 
First recall that we say a real valued random variable $X$ has
Logistic distribution if its distribution function is given by 
(\ref{Equ:Logistic-defi}), namely, 
\[
H(x) = \bP\left( X \leq x \right) = \frac{1}{1+e^{-x}}, \,\,\, x \in \Rbold.
\]
The following facts hold for the function $H$. 

\begin{Fact}
\label{fact:H-1}
$H$ is infinitely differentiable, and 
$H'(\cdot) = H(\cdot) \Hbar(\cdot)$, where
$\Hbar(\cdot) = 1 - H(\cdot)$.
\end{Fact}

\noindent
\emph{Proof :} From the definition it follows that $H$ is infinitely
differentiable on $\Rbold$. Further, 
\begin{eqnarray*}
H'(x) & = & \frac{1}{1+e^{-x}} \times \frac{e^{-x}}{1+e^{-x}} \\
      &   &  \\
      & = & H(x) \, \Hbar(x) \,\,\, \forall \,\,\, x \in \Rbold
\end{eqnarray*}
\qed

\begin{Fact}
\label{fact:H-2}
$H$ is symmetric around 0, that is, 
$H(-x) = \Hbar(x) \,\, \forall \,\, x \in \Rbold$.
\end{Fact}

\noindent
\emph{Proof :} From the definition we get that for any $x \in \Rbold$, 
\[
H(-x) = \frac{1}{1+e^x} = \frac{e^{-x}}{1+e^{-x}} = \Hbar(x).
\]
\qed

\begin{Fact}
\label{fact:H-3}
$\Hbar$ is the unique solution of the non-linear integral equation
\begin{equation} \tag{A1}
\Hbar(x) = \exp\left(- \int_{-x}^{\infty} \! \Hbar(s) \, ds 
\right), \,\,\, \forall \,\,\, x \in \Rbold.
\end{equation} 
\end{Fact}

\noindent
\emph{Proof :} Notice that the equation (A1) is nothing but Logistic RDE, 
this is because
\[
\bP\left( \min_{j \geq 1} \left( \xi_j - X_j \right) > x \right) 
= \exp\left(- \int_{-x}^{\infty} \! \Hbar(s) \, ds 
\right), \,\,\,\, \forall \,\,\, x \in \Rbold
\]
where $\left(X_j\right)_{j \geq 1}$ are i.i.d. with distribution function
$H$ and are independent of $\left(\xi_j\right)_{j \geq 1}$, which are
points of a Poisson point process of rate $1$ on $(0,\infty)$.
Thus from the fact that $\Hbar$ is the unique solution of Logistic RDE
(Lemma 5 of \cite{Al01})
we conclude that $\Hbar$ is unique solution of equation (A1). \qed

\section*{Acknowledgments}

This work was done in University of California, Berkeley, 
as a part of the author's doctoral dissertation, written under 
guidance of Professor David J. Aldous, whom the author would like to 
thank for suggesting the problem and for many illuminating discussion. 
The author would also like to thank an anonymous referee for 
some useful comments on an earlier version of the paper.

\bibliography{logistic.bib}

\begin{thebibliography}{10}

\bibitem{Al92a}
David Aldous.
\newblock Asymptotics in the random assignment problem.
\newblock {\em Probab. Theory Related Fields}, 93(4):507--534, 1992.

\bibitem{AlSt04}
David Aldous and J.~Michael Steele.
\newblock The objective method: probabilistic combinatorial optimization and
  local weak convergence.
\newblock In {\em Probability on discrete structures}, volume 110 of {\em
  Encyclopaedia Math. Sci.}, pages 1--72. Springer, Berlin, 2004.

\bibitem{Al00}
David~J. Aldous.
\newblock The percolation process on a tree where infinite clusters are frozen.
\newblock {\em Math. Proc. Cambridge Philos. Soc.}, 128(3):465--477, 2000.

\bibitem{Al01}
David~J. Aldous.
\newblock The $\zeta(2)$ {L}imit in the {R}andom {A}ssignment {P}roblem.
\newblock {\em Random Structures Algorithms}, 18(4):381--418, 2001.

\bibitem{AlBa05}
David~J. Aldous and Antar Bandyopadhyay.
\newblock A survey of max-type recursive distributional equations.
\newblock {\em Ann. Appl. Probab.}, 15(2):1047--1110, 2005.

\bibitem{Ba05a}
Antar Bandyopadhyay.
\newblock Hard-{C}ore {M}odel on {R}andom {G}raphs.
\newblock (preprint), 2005.

\bibitem{Ba05b}
Antar Bandyopadhyay.
\newblock A {N}ecessary and {S}ufficient {C}ondition for the
  {T}ail-{T}riviality of a {R}ecursive {T}ree {P}rocess.
\newblock (preprint, available at {\tt
  <http://www.arxiv.org/pdf/math.PR/0511203>}), 2005.

\bibitem{BaGa05}
Antar Bandyopadhyay and David Gamarnik.
\newblock Counting without sampling. {N}ew algorithms for enumeration problems
  using statistical physics.
\newblock To appear in the Proceedings of the ACM-SIAM Symposium on Discrete
  Algorithms 2006, (available at {\tt
  <http://www.arxiv.org/pdf/math.PR/0510471>}), 2006.

\bibitem{Bill95}
Patrick Billingsley.
\newblock {\em Probability and measure}.
\newblock John Wiley \& Sons Inc., New York, third edition, 1995.
\newblock A Wiley-Interscience Publication.

\bibitem{FiJa00}
James~Allen Fill and Svante Janson.
\newblock A characterization of the set of fixed points of the {Q}uicksort
  transformation.
\newblock {\em Electron. Comm. Probab.}, 5:77--84 (electronic), 2000.

\bibitem{Ga04}
D.~Gamarnik, T.~Nowicki, and G.~Swirscsz.
\newblock Maximum {W}eight {I}ndependent {S}ets and {M}atchings in {S}parse
  {R}andom {G}raphs. {E}xact {R}esults using the {L}ocal {W}eak {C}onvergence
  {M}ethod.
\newblock To appear in \emph{Random Structures and Algorithms}, (available at
  {\tt <http://www.arxiv.org/pdf/math.PR/0309441>}), 2004.

\bibitem{LW03}
Svante Linusson and Johan W{\"a}stlund.
\newblock A proof of {P}arisi's conjecture on the random assignment problem.
\newblock {\em Probab. Theory Related Fields}, 128(3):419--440, 2004.

\bibitem{MePa87}
M.~M\'{e}zard and G.~Parisi.
\newblock On the solution of the random link matching problem.
\newblock {\em J. Physique}, 48:1451--1459, 1987.

\bibitem{MePaVi87}
M.~M\'{e}zard, G.~Parisi, and M.~A. Virasoro.
\newblock Spin {G}lass {T}heory and {B}eyond.
\newblock {\em World Scientific, Singapore}, 1987.

\bibitem{NPS03}
C.~Nair, B.~Prabhakar, and M.~Sharma.
\newblock A proof of the conjecture due to {P}arisi for the finite random
  assignment problem.
\newblock (preprint, available at {\tt
  http://www.stanford.edu/$\sim$balaji/papers/parisi.pdf}), 2003.

\bibitem{Pa98}
Giorgio Parisi.
\newblock A conjecture on random bipartite matching.
\newblock {\em Physics e-Print archive, {\tt
  http://xxx.lanl.gov/ps/cond-mat/9801176}}, 1998.

\bibitem{RosRu01}
U.~R{\"o}sler and L.~R{\"u}schendorf.
\newblock The contraction method for recursive algorithms.
\newblock {\em Algorithmica}, 29(1-2):3--33, 2001.
\newblock Average-case analysis of algorithms (Princeton, NJ, 1998).

\bibitem{Ros92}
Uwe R{\"o}sler.
\newblock A fixed point theorem for distributions.
\newblock {\em Stochastic Process. Appl.}, 42(2):195--214, 1992.

\bibitem{Sim72}
George~F. Simmons.
\newblock {\em Differential equations with applications and historical notes}.
\newblock McGraw-Hill Book Co., New York, 1972.
\newblock International Series in Pure and Applied Mathematics.

\bibitem{St97}
J.~Michael Steele.
\newblock {\em Probability theory and combinatorial optimization}.
\newblock Society for Industrial and Applied Mathematics (SIAM), Philadelphia,
  PA, 1997.

\end{thebibliography}

\end{document}